\documentclass[11pt]{article}

\usepackage{a4}
\usepackage{fullpage}
\usepackage{amsfonts}
\usepackage{amssymb}
\usepackage{amsmath}
\usepackage{amsthm}
\usepackage{amsopn}
\usepackage{graphicx}
\usepackage{times}

\newcommand{\finishline}{\hfill\hbox{}\linebreak[4]}

\newcommand{\Z}{\mathbb Z}

\newcommand{\R}{\mathbb R}

\newcommand{\GL}{\mathsf{GL}}
\newcommand{\BA}{{\mathfrak A}}
\newcommand{\BB}{{\mathfrak B}}
\newcommand{\BE}{{\mathfrak E}}

\newcommand{\MC}{{\cal C}}

\newcommand{\MF}{{\cal F}}
\newcommand{\MM}{{\cal M}}

\newcommand{\MV}{{\cal V}}
\renewcommand{\vec}[1]{\boldsymbol{#1}}

\DeclareMathOperator{\conv}{conv}

\DeclareMathOperator{\id}{id}

\DeclareMathOperator{\diag}{diag}

\DeclareMathOperator{\sgn}{sgn}
\DeclareMathOperator{\DV}{DV}

\theoremstyle{definition}
\newtheorem{definition}{Definition}[section]

\newtheorem{lemma}[definition]{Lemma}

\newtheorem{theorem}[definition]{Theorem}
\newtheorem{corollary}[definition]{Corollary}
\newtheorem{conjecture}[definition]{Conjecture}
\newtheorem{example}[definition]{Example}

\newtheorem{remark}[definition]{Remark}

\author{
Frank Vallentin
\footnote{Partially supported by the Edmund Landau Center for Research in Mathematical Analysis and Related Areas, sponsored by the Minerva Foundation (Germany).}\\[0.5ex]
}
\title{A NOTE ON SPACE TILING ZONOTOPES}

\begin{document}
 
\maketitle


\begin{abstract}
In 1908 Voronoi conjectured that every convex polytope which tiles space face-to-face by translations is affinely equivalent to the Dirichlet-Voronoi polytope of some lattice. In 1999 Erdahl proved this conjecture for the special case of zonotopes. A zonotope is a projection of a regular cube under some affine transformation. In 1975 McMullen showed several equivalent conditions for a zonotope to be a space tiling zonotope, i.e.\ a zonotope which admits a face-to-face tiling of space by translations. Implicitly, he
related space tiling zonotopes to a special class of oriented matroids (regular matroids). We will extend his result to give a new proof of Voronoi's conjecture for zonotopes using oriented matroids. This enables us to distinguish between combinatorial and metrical properties and
to apply the fact that oriented matroids considered here have an essentially unique realization. Originally, this is a theorem due to Brylawski and Lucas. By using oriented matroid duality we interpret a part of McMullen's arguments as an elegant geometric proof of this theorem in the special case of real numbers. This note grew apart \cite{Vallentin-00:Paralleloeder}.
\end{abstract}


\section{Introduction}

Let $V$ be a $d$-dimensional real vector space.
A \textit{parallelohedron} $P \subseteq V$ is a convex polytope which admits a face-to-face tiling of~$V$ by translations.
One construction to obtain parallelohedra is the following:
Let $L \subseteq V$ be a lattice (By a lattice $L \subseteq V$ we mean the $\Z$-span of $n \leq d$ linearly independent vectors.) and let $(\cdot, \cdot) : V \times V \to \R$ be an inner product, then the polytope
$$
\DV(L, (\cdot, \cdot)) =
\{\vec{x} \in V :
\mbox{for all $\vec{v} \in L$ we have $d(\vec{x}, \vec{0}) \leq
d(\vec{x}, \vec{v})$}\},\;\; d(\vec{x}, \vec{y}) = \sqrt{(\vec{x} - \vec{y}, \vec{x} - \vec{y})},
$$
is a parallelohedron. It is called \textit{Dirichlet-Voronoi polytope} of $(L, (\cdot, \cdot))$. The translates $\DV(L, (\cdot, \cdot)) + \vec{v}$, $\vec{v} \in L$, give a face-to-face tiling of~$V$. 

\begin{conjecture} (Voronoi's Conjecture, \cite{Voronoi-08:Recherches} p.\ $210$f)
\finishline
For every parallelohedron $P \subseteq V$ there exists a lattice $L
\subseteq V$ and an inner product $(\cdot, \cdot): V \times V \to \R$
such that the Dirichlet-Voronoi polytope of $(L, (\cdot, \cdot))$ is a
translate of~$P$.
\end{conjecture}

Voronoi's conjecture has been proven in several special cases although it is open in general. Quite recently, Erdahl \cite{Erdahl-99:Zonotopes} proved the conjecture for space tiling zonotopes. 

In this paper we will give a new proof of Voronoi's conjecture for zonotopes. We will extend the paper \cite{McMullen-75:Space} of McMullen in which he gave several equivalent conditions for a zonotope to be space tiling. We mainly add the concept of zonotopal lattices and a somewhat fiddly matrix computation to McMullen's paper. The main difference of our proof to Erdahl's is that we are using oriented matroid terminology (which is more than only implicit in McMullen's paper and the connection between space tiling zonotopes and regular matroids has been observed by several authors). This enables us to distinguish between combinatorial and metrical properties and to make use of the fact that the matroids we are considering have a unique realization up to projective transformations (see Remark~\ref{rem:unique}). Furthermore, we do not need to switch to dual tilings, ``lattice dicings'' in the terminology of Erdahl. Deza and Grishukhin gave another proof of Voronoi's conjecture for zonotopes in \cite{Deza-Grishukhin-02:Once} where they actually used oriented matroids. On the other hand they had to apply a characterization of parallelohedra of Venkov and McMullen whose proof is rather elaborate which we will not need. 


\section{Notation}

In this section we collect the notations we will use. For more information see e.g.\ \cite{Bjoerner-93:Oriented} and \cite{Ziegler-95:Lectures}. 

\paragraph{Zonotopes and Oriented Matroids.}

Let $X = (\vec{x}_1, \ldots, \vec{x}_n) \in V^n$ be a vector configuration spanning~$V$. The vector configuration defines a zonotope $Z(X) = S_1 + \cdots + S_n$ with $S_i = \conv\{\vec{x}_i, -\vec{x}_i\}$, $i = 1, \ldots, n$.
It also gives \textit{vectors} 
$\MV(X) = \{(\sgn \lambda_1, \ldots, \sgn \lambda_n) \in \{-1,0,+1\}^n: \sum_{i=1}^d \lambda_i \vec{x}_i = \vec{0}\}$ and \textit{covectors} $\MV^*(X) = \{(\sgn f(\vec{x}_1), \ldots, \sgn f(\vec{x}_n)) \in \{-1,0,+1\}^n: \mbox{$f : V \to \R$ linear functional}\}$ of the oriented matroid $\MM(X)$.
By $\MF$ we denote the set of all $k$-dimensional faces, $k = 0, \ldots, d$, of the zonotope $Z(X)$. The map 
$$
\varphi : \MV^*(X) \to \MF,\; 
(v_1, \ldots, v_n) \mapsto
\sum_{i : v_i > 0} \vec{x}_i - \sum_{i : v_i < 0} \vec{x}_i + \sum_{i : v_i = 0} S_i
$$
defines an order-reversing bijection between $\MF$ ordered by inclusion and $\MV^*(X) \subseteq \{-1,0,+1\}^n$ ordered by componentwise extension of $0 \prec \pm 1$. The minimal elements in the poset $(\MV^*(X), \prec)$ are called \textit{cocircuits} (Notation: $\MC^*(X)$) of the oriented matroid.
Let $B \subseteq X$ be a subset of $X$ that forms a basis of $V$ and let $\vec{x}_j \in B$ some basis vector. Then there exist exactly two cocircuits $\pm (c_1, \ldots, c_n) \in \MC^*(X)$ with $c_j \neq 0$ and $c_i = 0$ if $i \in B \backslash \{\vec{x}_j\}$. They are called \textit{basis cocircuits} of~$j$ with respect to~$B$ (Notation: $c^*(j,B)$).
Two vector configurations $X = (\vec{x}_1, \ldots, \vec{x}_n)$ and $Y = (\vec{y}_1, \ldots, \vec{y}_n)$ are called \textit{projectively equivalent} if there exists a regular transformation $A \in \GL(V)$ and a scaling transformation $D$ ($\vec{x}_i D = \lambda_i\vec{x}_i$, $\lambda_i \neq 0$) with $AXD = Y$. In this situation we have $\MM(X) = \MM(Y)$. 
For notational convenience we assume that the vector space $V$ is $\R^d$ and we interpret $X$ as a matrix with $d$ rows an $n$ columns. With $X$ we associate a vector configuration $\overline{X} = (\overline{\vec{x}}_1, \ldots, \overline{\vec{x}}_n) \in \R^{(n-d) \times n}$ of rank~$n-d$ so that the matrix $\left(\begin{smallmatrix}X\\\overline{X}\end{smallmatrix}\right) \in \R^{n \times n}$ has rank $n$ and the first $n$ rows are orthogonal to the last $n-d$ rows. In this case we have $\MV^*(\overline{X}) = \MV(X)$ and $\MV^*(X) = \MV(\overline{X})$. Note that $\overline{X}$ is unique up to a linear transformation.

\paragraph{Lattices.}

A lattice vector $\vec{v} \in L \backslash \{\vec{0}\}$ is called \textit{strict Voronoi vector} if $\pm \vec{v}$ are the only vectors in the coset $\vec{v} + 2L$ with minimal norm. By a theorem of Voronoi (see e.g.\ \cite{Voronoi-08:Recherches}), we have 
$$
\DV(L, (\cdot, \cdot)) = \bigcap\limits_{\mbox{$\vec{v}$ strict Voronoi vector}} \{\vec{x} \in V : (\vec{x}, \vec{v}) \leq \frac{1}{2}(\vec{v},\vec{v})\},
$$
and any of these inequalities defines a facet of $\DV(L, (\cdot, \cdot))$.
A lattice $L \subseteq \Z^n$ gives the covectors of the oriented matroid $\MM(L)$ by $\MV^*(L) = \{(\sgn v_1, \ldots, \sgn v_n) \in \{-1,0,+1\}^n: \vec{v} \in L\}$. The \textit{support} of
$\vec{v} \in L$ is $\underline{\vec{v}} = \{i \in
\{1,\ldots,n\} : v_i \neq 0\}$. The vector $\vec{v}$ is called
\textit{elementary} if $\vec{v} \in \{-1,0,+1\}^n \backslash
\{\vec{0}\}$ and if $\vec{v}$ has minimal support among all vectors in~$L \backslash \{\vec{0}\}$. We say that two 
vectors $\vec{v}, \vec{w} \in L$ are \textit{conformal} if $v_i
\cdot w_i \geq 0$ for all $i \in \{1, \ldots, n\}$.

\paragraph{Linear Algebra.}

Let $W$ be another real vector space of finite dimension, and let $f : V \to W$ be a linear map. Recall that $M^\BA_\BB(f)$ denotes the matrix representing~$f$ with respect to the basis~$\BA$ of~$V$ and the basis~$\BB$ of~$W$.


\section{Zonotopal Lattices}
\label{sec:zl}

\begin{definition}
Let $L \subseteq \Z^n$ be a lattice, and let $(\cdot,
\cdot)$ be an inner product on $\R^n$ where the canonical basis $\BE = (\vec{e}_1, \ldots, \vec{e}_n)$ forms an orthogonal basis (but not necessarily an orthonormal basis). The pair $(L, (\cdot, \cdot))$
is called \textit{zonotopal lattice} if for every vector $\vec{v} \in L \backslash \{\vec{0}\}$ there exists an elementary vector $\vec{u} \in L$ with $\underline{\vec{u}} \subseteq \underline{\vec{v}}$.
\end{definition}

From now on until the end of Section~\ref{sec:zl} we assume that $(L,
(\cdot, \cdot))$ is a zonotopal lattice and that $F$ is the
subspace spanned by $L$. We will need the following elementary properties of zonotopal lattices. Proofs of these facts can be found e.g.\ in
\cite{Tutte-71:Introduction}, Chapter 1.2.

\begin{lemma}
The elementary vectors of $L$ are cocircuits of $\MM(L)$ and vice versa.
Every vector of $L$ can be written as a sum of pairwise conformal elementary vectors.
If $\vec{v}$ is an elementary vector, and we have $\vec{u} \in L$ with
$\underline{\vec{u}} = \underline{\vec{v}}$, then there exists a
factor $\alpha \in \Z \backslash \{-1,+1\}$ such that $\vec{u} =
\alpha \vec{v}$.
\end{lemma}

The aim of the rest of the section is to show that Dirichlet-Voronoi
polytopes of zonotopal lattices are indeed zonotopes. Furthermore we
will show that the oriented matroid $\MM(L)$ and the oriented matroid
given by the vector configuration of the zonotope $\DV(L, (\cdot, \cdot))$
coincide. These facts are stated in \cite{Loesch-90:Reduktionstheorie}
without proofs.

\begin{lemma}
\label{lem:cocircuits}
Every elementary vector of~$L$ is a strict Voronoi vector of $L$ and
vice versa.
\end{lemma}

\begin{proof}
Let $\vec{v} \in L$ be a strict Voronoi vector. Write $\vec{v} =
\sum_{i=1}^m \vec{w}_i$ as a sum of pairwise conformal elementary vectors
$\vec{w}_i \in L$.  Assume that $m \geq 2$. Defining $\vec{u} =
\vec{v} - 2\vec{w}_1$ gives $\vec{u} \neq \pm\vec{v}$ and $(\vec{u},
\vec{u}) = (\vec{v}, \vec{v}) - 4(\vec{v} - \vec{w}_1, \vec{w}_1)$.
Since the vectors $\vec{w}_i$, $i = 1, \ldots, m$, are pairwise conformal we have $(\vec{v} - \vec{w}_1,
\vec{w}_1) \geq 0$, and $\pm\vec{v}$ is not the unique shortest vector
in $\vec{v} + 2L$. In this case $\vec{v}$ cannot be a
strict Voronoi vector. Hence, $m = 1$.

Let $\vec{v} \in L$ be an elementary vector, and let $\vec{u} \in
\vec{v} + 2L$ be a lattice vector with $\vec{u} \neq \pm \vec{v}$.
We have $\vec{v} - \vec{u} \in 2L \subseteq 2\Z^d$ and $\vec{v}_i \in
\{-1,0,+1\}$, which shows $\underline{\vec{v}} \subseteq \underline{\vec{u}}$.
The case $\underline{\vec{v}} \neq \underline{\vec{u}}$ immediately
leads to $(\vec{v}, \vec{v}) < (\vec{u}, \vec{u})$.  If
$\underline{\vec{v}} = \underline{\vec{u}}$, then there exists a factor
$\alpha \in \Z \backslash \{-1,+1\}$ so that $\vec{u} = \alpha
\vec{v}$, hence $(\vec{v}, \vec{v}) < (\vec{u}, \vec{u})$.  In both
cases $\pm\vec{v}$ are the only shortest vectors in $\vec{v} + 2L$.
\end{proof}

\begin{lemma} (see e.g.\ Theorem 22.6 in \cite{Rockafellar-70:Convex})
\label{lem:farkas}
\finishline
Let $\vec{x} \in \R^n$ be a vector, and let
$\alpha_1, \ldots, \alpha_n \in \R \cup \{\pm\infty\}$ be real
numbers.  Either there exists a vector $\vec{y}' \in F^{\perp}$ lying in the parallelotope
$\vec{x} + \prod_{i=1}^n [-\alpha_i, \alpha_i]$,
or there exists a vector $\vec{y} \in F$ such that for all $\vec{z} \in
\vec{x} + \prod_{i=1}^n [-\alpha_i, \alpha_i]$ the inequality $(\vec{y}, \vec{z}) > 0$ holds.
If the second condition holds, then one can choose $\vec{y}$ to be an
elementary vector of $L$.
\end{lemma}

\begin{theorem}
\label{th:zonotopal}
By $\pi : \R^n \to F$ denote the orthogonal projection of $\R^n$ onto
$F$. Then, $\DV(L) = \pi(\DV(\Z^n)) =
\pi([-\frac{1}{2}, +\frac{1}{2}]^n)$.
\end{theorem}

\begin{proof}
Suppose that $\vec{x} \in \DV(\Z^n)$. For all $\vec{v} \in \Z^n
\backslash \{\vec{0}\}$ the inequality $(\vec{x}, \vec{v}) \leq 
\frac{1}{2}(\vec{x}, \vec{v})$ holds. Write $\vec{x} = \vec{y} + \vec{y}'$ with $\vec{y} = \pi(\vec{x}) \in F$ and
$\vec{y}' \in F^{\perp}$. For all $\vec{v} \in L \backslash
\{\vec{0}\}$ we have $(\vec{y}, \vec{v}) = (\vec{x}, \vec{v}) - (\vec{y}', \vec{v}) \leq
\frac{1}{2}(\vec{v}, \vec{v})$. Thus, $\pi(\vec{x}) \in \DV(L)$.

Suppose now that $\vec{y} \in \DV(L)$. If there exists $\vec{x} \in
(-\vec{y} + \prod_{i=1}^n [-\frac{1}{2}, \frac{1}{2}]) \cap
F^{\perp}$, then $\vec{y} + \vec{x} \in
\DV(\Z^n)$ and $\pi(\vec{y} + \vec{x}) =
\vec{y}$. Assume that such a vector does not exist. Then by Lemma~\ref{lem:farkas} 
there is an elementary lattice vector $\vec{v} \in L$ so that
$(\vec{v}, -\vec{y} + \prod_{i=1}^n [-\frac{1}{2}, \frac{1}{2}]) >
0$. This implies $(\vec{v}, -\vec{y} - \frac{1}{2} \vec{v}) > 0$,
respectively $(\vec{v}, \vec{y}) < \frac{1}{2}(\vec{v}, \vec{v})$, in
contradiction to the assumption $\vec{y} \in \DV(L)$
\end{proof}

\begin{example}
Figure~$1$ shows the situation for the zonotopal lattice $(L = \Z(1,1,0) + \Z(0,1,1), \langle \cdot, \cdot \rangle)$ where $\langle \cdot, \cdot \rangle$ denotes the standard inner product. The light grey plane is~$F$ and the dark grey hexagon is $\DV(L, \langle \cdot, \cdot \rangle)$. 
\begin{center}
\begin{center}
\unitlength1cm
\begin{picture}(10,8.8)
\put(0,0){\includegraphics[width = 10cm]{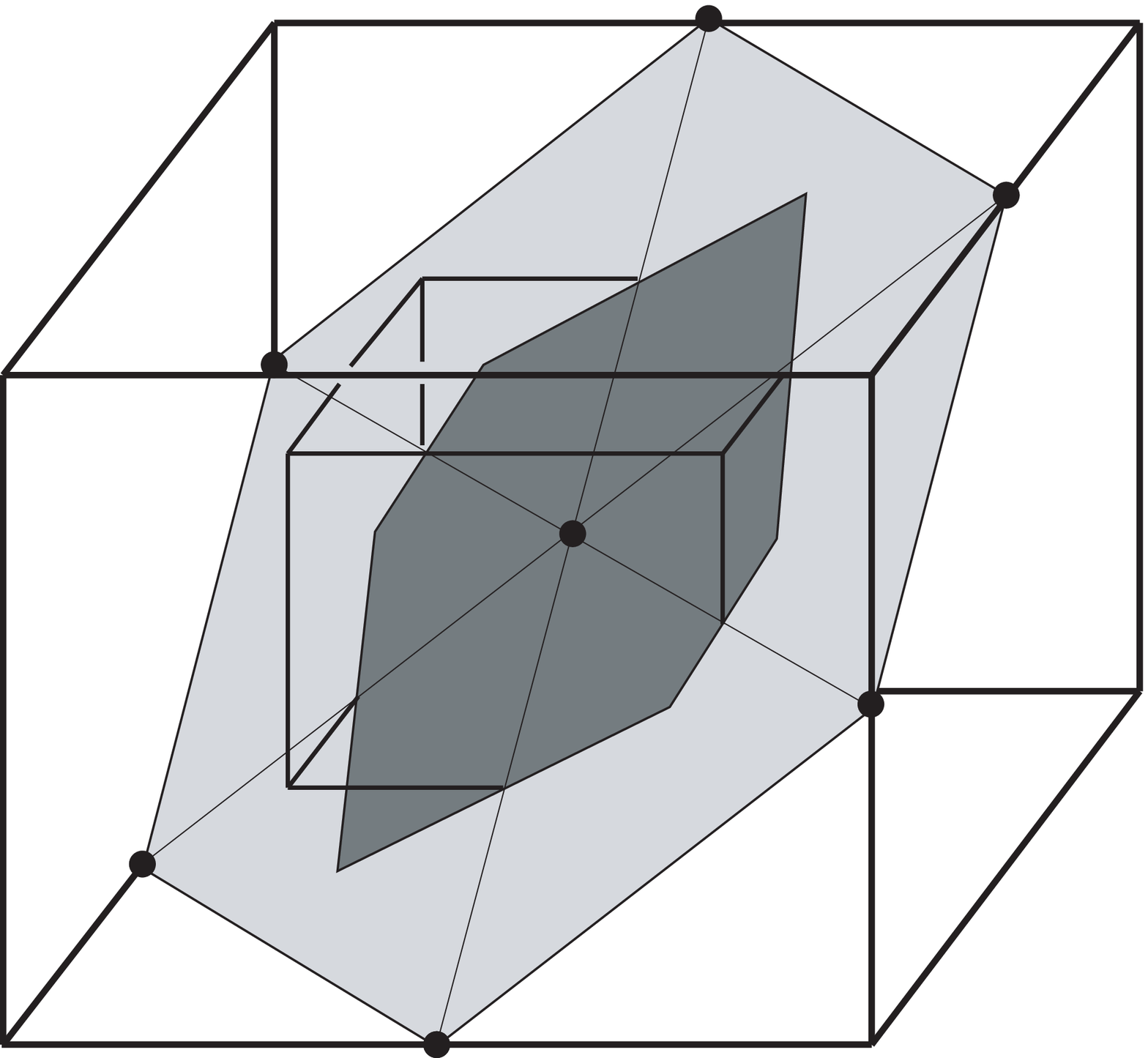}}
\put(4.7,4){\footnotesize $(0,0,0)$}
\put(5.4,8){\footnotesize $(-1,0,1)$}
\put(7.5,3.6){\footnotesize $(1,1,0)$}
\put(1.4,6.2){\footnotesize $(-1,-1,0)$}
\put(1.1,2.2){\footnotesize $(0,-1,-1)$}
\put(3.5,0.7){\footnotesize $(1,0,-1)$}
\put(7.7,7.5){\footnotesize $(0,1,1)$}
\end{picture}
\end{center}
\smallskip
\textbf{\textsf{Figure 1.}} Dirichlet-Voronoi Polytope of a Zonotopal Lattice.
\end{center} 
\end{example}

From Theorem~\ref{th:zonotopal} and Lemma~\ref{lem:cocircuits} it follows
that the cocircuits of $\MM(L)$ and the oriented
matroid which belongs to the vector configuration $Y = (\pi(\vec{e}_1),
\ldots, \pi(\vec{e}_n)) \in
\R^{n \times n}$ defining the zonotope~$\DV(L)$ coincide:

\begin{corollary}
The two oriented matroids $\MM(L)$ and $\MM(Y)$ coincide.
\end{corollary}


\section{Equivalent Conditions for Space Tiling Zonotopes}

In this section we give a new proof of Voronoi's conjecture for zonotopes which was first proven by Erdahl \cite{Erdahl-99:Zonotopes}.
McMullen \cite{McMullen-75:Space}, extending works of Coxeter \cite{Coxeter-62:Classification} and Shephard \cite{Shephard-74:Space}, proved several equivalent conditions for space tiling zonotopes. He showed that the vector configuration associated to a space tiling zonotope is projectively equivalent to the vector configuration of the Dirichlet-Voronoi polytope of a lattice. But he did not show Voronoi's conjecture for zonotopes. Here we extend McMullen's proof by a few arguments. This yields a proof of Voronoi's conjecture for zonotopes. 
For organizing our proof we use McMullen's original numbering with minor changes. We interpret condition (III) in our terminology and augment it by a statement which is given in \cite{McMullen-75:Space} without number. Our condition (VII) is a strengthening of McMullen's original (VII), a reformulation of Voronoi's conjecture for zonotopes.

We give an outline of the proof first. Starting from a space tiling zonotope $Z(X)$ we define the lattice $L$ spanned by the cocircuits $\MC^*(X)$.
We show that $L$ together with any suitable inner product $\langle \cdot, \cdot \rangle$ is a zonotopal lattice. Furthermore, the oriented matroids $\MM(X)$ and $\MM(L)$ coincide. Since these oriented matroids have an essentially unique realization, $X$ and the vector configuration associated to $\DV(L, \langle \cdot, \cdot \rangle)$ are projectively equivalent. Finally, by a matrix computation we show that we can find an inner product $(\cdot, \cdot)$ so that $\DV(L, (\cdot, \cdot))$ is affinely equivalent to $Z(X)$, i.e.\ there is an affine map $f$ so that $f(Z(X)) = \DV(L, (\cdot, \cdot))$.

\begin{theorem}
Let $X = (\vec{x}_1, \ldots, \vec{x}_n) \in \R^{d \times n}$ be a $d$-dimensional vector configuration.
By $Z = Z(X)$ we denote the zonotope defined by~$X$.
Let $L = \langle \MC^* \rangle_\Z \subseteq
\Z^n$ be the lattice spanned by the cocircuits of~$\MM(X)$. Let $F$ be the subspace spanned by $L$, and let $\pi : \R^n
\to F$ be the orthogonal projection with respect to the standard inner product $\langle \cdot, \cdot \rangle$. Define the vector configuration
$Y = (\vec{y}_1, \ldots, \vec{y}_n)$ by $\vec{y}_i = \pi(\vec{e}_i)$,
$i = 1, \ldots, n$.  Then, the following statements are equivalent.
\begin{enumerate}
\item[(I)] $Z$ tiles $\R^d$ face-to-face by translations.
\item[(II)] Every $(d-2)$-space spanned by vectors of $X$ is contained in $2$ or $3$ hyperplanes spanned by vectors of $X$. 
\item[(III)] The lattice $(L, \langle \cdot, \cdot\rangle)$ is zonotopal of rank $d$. The vector configurations $X$ and $Y$ are projectively equivalent.
\item[(VII)] There exists an inner product $(\cdot, \cdot)$ so that $\DV(L, (\cdot, \cdot))$ and $Z$ are affinely equivalent.
\end{enumerate}
\end{theorem}

\begin{proof}
\finishline
$\underline{\mbox{(I)} \Longrightarrow \mbox{(II)}:}\quad$ 
See \cite{Coxeter-62:Classification}, $\S 9$.
This corresponds to the fact that two-dimensional space tiling zonotopes are either quadrangles or centrally symmetric hexagons.

\noindent $\underline{\mbox{(II)} \Longrightarrow \mbox{(III)}:}\quad$ In our notions, Lemma 4 of \cite{McMullen-75:Space} translates into the fact that~$L$ has rank~$d$. Thus, the oriented matroids $\MM(L)$ and
$\MM(X)$ have the same rank. The covectors of $\MM(X)$ are
covectors of $\MM(L)$ by definition. From a rank argument we get $\MM(X) = \MM(L)$. The cocircuits are the minimal elements in the poset of covectors. Cocircuits of~$\MM(L)$ correspond to elementary vectors and covector of~$\MM(L)$ correspond to non-trival vectors of~$L$. Thus, the
support of every vector of $L \backslash \{\vec{0}\}$ is contained in
the support of some elementary vector of $L$. This implies that $(L, \langle \cdot, \cdot \rangle$ is a zonotopal lattice.

Suppose (after suitable relabelling, if necessary) that $B =
(\vec{x}_1, \ldots, \vec{x}_d)$ is a basis of $\R^d$.  For $j \in \{1,
\ldots, d\}$ consider the basic cocircuit $c^*(j, B)$ of $j$ with
respect to $B$ and with $c^*(j, B)_j = +1$. There are $\lambda_{ij}$
with $\overline{\vec{x}}_j + \sum_{i=d+1}^n
\lambda_{ij} \overline{\vec{x}}_i = \vec{0}$. Define $\overline{\vec{x}}'_j = \overline{\vec{x}}_j$, $j = d+1,
\ldots, n$ and $\overline{\vec{x}}'_j = -\sum_{i=d+1}^{n} \sgn\lambda_{ij} \overline{\vec{x}}_i$, $j = 1, \ldots, d$.
Note that $\sgn\lambda_{ij} = c^*(j, B)_i$.  On p.\ $209$ of \cite{McMullen-75:Space}
McMullen proves that the vector configurations
$(\overline{\vec{x}}_1, \ldots,
\overline{\vec{x}}_n)$ and $(\overline{\vec{x}}'_1, \ldots,
\overline{\vec{x}}'_n)$ are projectively equivalent. 
In the same way we construct the vectors $\overline{\vec{y}}'_i$ from
$\overline{\vec{y}}_i$, $i = 1, \ldots, n$, and 
the same argument shows that the vector configurations
$(\overline{\vec{y}}_1, \ldots,
\overline{\vec{y}}_n)$ and $(\overline{\vec{y}}'_1, \ldots,
\overline{\vec{y}}'_n)$ are projectively equivalent.
Since the cocircuits of $\MM(X)$ and $\MM(Y)$ coincide and by
construction of the vectors $\overline{\vec{x}}'_i$ and
$\overline{\vec{y}}'_i$, $i = 1, \ldots, d$, we see that
$(\overline{\vec{x}}'_1, \ldots, \overline{\vec{x}}'_n)$ and
$(\overline{\vec{y}}'_1, \ldots,
\overline{\vec{y}}'_n)$ are projectively equivalent. By
transitivity the vector configurations $(\overline{\vec{x}}_1, \ldots,
\overline{\vec{x}}_n)$ and $(\overline{\vec{y}}_1, \ldots,
\overline{\vec{y}}_n)$ are projectively equivalent, too.
Then, by duality, the vector configurations $X$ and $Y$ are
projectively equivalent.

\noindent $\underline{\mbox{(III)} \Longrightarrow \mbox{(VII)}:}\quad$ 
We extend $X \in
\R^{d \times n}$ to a square matrix with $n$~rows by adding rows with zeroes.
By assumption we have the equation $X = AYD$ for a suitable matrix $A \in \GL_n(\R)$
and a diagonal matrix $D = \diag(\lambda_1, \ldots, \lambda_n)$ with
strictly positive diagonal elements. Choose the inner product
on $\R^n$ to be $(\vec{e}_i, \vec{e}_j) = \lambda_i \delta_{ij}$ where $\delta_{ij}$ is the
Kronecker delta. With $\pi'$ denote the orthogonal projection of
$\R^n$ onto $F$ with respect to this inner product. In the next
paragraph we will construct a matrix $B \in
\GL_n(\R)$ so that $BYD = (\pi'(\vec{e}_1), \ldots, \pi'(\vec{e}_n))$. Then,
$X = A B^{-1} (\pi'(\vec{e}_1), \ldots, \pi'(\vec{e}_n))$, which means
that the zonotopes $Z$ and $\DV(L, (\cdot, \cdot))$ are affinely
equivalent.

Now let us construct $B$. Let $\BA$ be a basis of the subspace $F$,
and let $\BA'$ be a basis of the subspace orthogonal to $F$ with
respect to the standard inner product $\langle \cdot, \cdot
\rangle$, i.e.\ $\langle \vec{e}_i, \vec{e}_j \rangle = \delta_{ij}$. Then, $D^{-1}\BA'$ is a basis of the subspace orthogonal
to~$F$ with respect to $(\cdot, \cdot)$. We have
$$
X = M^\BE_\BE(\pi) = M^{(\BA, \BA')}_\BE(\id) \diag(1, \ldots, 1, 0, \ldots, 0)
\left(M^{(\BA, \BA')}_\BE(\id)\right)^{-1},
$$
and
$$
(\pi'(\vec{e}_1), \ldots, \pi'(\vec{e}_n)) = M^\BE_\BE(\pi') = M^{(\BA, D^{-1}\BA')}_\BE(\id) 
\diag(1, \ldots, 1, 0, \ldots, 0)
\left(M^{(\BA, D^{-1}\BA')}_\BE(\id)\right)^{-1}.
$$
Define 
$$
C =
(M^{(\BA,\BA')}_\BE(\id))^{-1} D M^{(\BA,D^{-1}\BA')}_\BE(\id). 
$$
Then, $\diag(1,\ldots, 1,0,\ldots,0) C$ is a block matrix, where only
the main minor of size $d \times d$, which has rank $d$, has non-zero entries.
Set
$$
C_1 =
\left({\diag(1,\ldots, 1,0,\ldots,0) C + \diag(0, \ldots, 0, 1, \ldots,
1)}\right)^{-1},
$$
and
$$ 
B = M^{(\BA,D^{-1}\BA')}_\BE(\id) C_1
\left(M^{(\BA,\BA')}_\BE(\id)\right)^{-1}.
$$
Hence,
\begin{eqnarray*}
B M^\BE_\BE(\pi) D & = & M^{(\BA,D^{-1}\BA')}_\BE(\id) C_1
\left(M^{(\BA,\BA')}_\BE(\id)\right)^{-1} 
M^\BE_\BE(\pi) D\\
& = & 
M^{(\BA,D^{-1}\BA')}_\BE(\id) C_1
\diag(1, \ldots, 1, 0, \ldots, 0)
\left(M^{(\BA, \BA')}_\BE(\id)\right)^{-1} D\\
& = & 
M^{(\BA,D^{-1}\BA')}_\BE(\id) C_1
\diag(1, \ldots, 1, 0, \ldots, 0) 
C \left(M^{(\BA, D^{-1}\BA')}_\BE(\id)\right)^{-1}\\
& = & 
M^{(\BA,D^{-1}\BA')}_\BE(\id) 
\diag(1, \ldots, 1, 0, \ldots, 0) 
\left(M^{(\BA, D^{-1}\BA')}_\BE(\id)\right)^{-1} \\
& = & M^\BE_\BE(\pi').
\end{eqnarray*}

\noindent $\underline{\mbox{(VII)} \Longrightarrow \mbox{(I)}:}\quad$
This assertion is trivial.
\end{proof}

\begin{remark}
\label{rem:unique}
The oriented matroid $\MM(L)$ of a zonotopal lattice $(L, (\cdot, \cdot))$ is a regular matroid (a comprehensive reference for regular matroids is \cite{Tutte-71:Introduction}). On the other hand for every regular matroid $\MM$ one can define a zonotopal lattice $(L, (\cdot, \cdot))$ with $\MM = \MM(L)$.
Regular matroids are exactly those matroids having a realization over every field. A theorem of Brylawski and Lucas \cite{Brylawski-Lucas-76:Uniquely} says that one can transform any two vector configurations realizing the same regular matroid into each other by a projective transformation. This fact is highly wrong for general matroids. The arguments in ``(II) $\Longrightarrow$ (III)'' give a geometric proof of this theorem in the special case of real numbers.
\end{remark}



\vskip 20 pt

\noindent
\textit{Einstein Institute of Mathematics, The Hebrew University of Jerusalem, Israel}

\noindent
\textit{e-mail:} \texttt{vallenti@ma.tum.de}

\end{document}